\documentclass[12pt]{article}
\usepackage{amssymb,amsmath, amsthm}
\newtheorem{thm}{Theorem}
\newtheorem*{theo}{Theorem}

\newtheorem{lemma}{Lemma}
\newenvironment{defin}{\medskip\noindent{\sc
Definition}.}{\goodbreak\medskip}
\newenvironment{nota}{\medskip\noindent{\sc
Notations}.}{\goodbreak\medskip}
\newenvironment{remk}{\noindent{\sc
Remark}.}{\goodbreak\vskip10pt}
\newtheorem{prop}{Proposition}

\def\demo{\medskip\goodbreak\noindent
     \hbox{\sc Proof \kern .3em}\ignorespaces}%
  \def \qedbox{$\square$}%
  \def \qed{\hglue1mm\hfill{\ifmmode\qedbox
     \else\unskip\ \hglue0mm\hfill\qedbox\medskip
      \goodbreak\fi}}%
\def\enddemo{\qed\goodbreak\vskip10pt}%

\def\qed{\hglue1mm\hfill\raise -2pt\hbox{\vrule\vbox to 10pt{\hrule width
4pt
                  \vfill\hrule}\vrule}}

\newcommand{\T}{\mathbb {T}}
\newcommand{\A}{\mathbb {A}}

\newcommand{\R}{\mathbb {R}}
\newcommand{\Q}{\mathbb {Q}}

\newcommand{\Z}{\mathbb {Z}}
\newcommand{\N}{\mathbb {N}}

\begin{document}
\title{
A non-differentiable essential  irrational invariant curve for a   $C^1$ symplectic twist map }
\author{M.-C. ARNAUD
\thanks{ANR KAM faible ANR-07-BLAN-0361}
\thanks{ANR DynNonHyp ANR BLAN08-2-313375}
\thanks{Universit\'e d'Avignon et des Pays de Vaucluse, Laboratoire d'Analyse non lin\' eaire et G\' eom\' etrie (EA 2151),  F-84 018Avignon,
France. e-mail: Marie-Claude.Arnaud@univ-avignon.fr}
}
\maketitle
\abstract{   We construct a $C^1$ symplectic twist map $f$ of the annulus that has an essential invariant curve $\Gamma$ such that:
\begin{enumerate}
\item[$\bullet$] $\Gamma$ is not differentiable;
\item[$\bullet$] the dynamic of $f_{|\Gamma}$ is conjugated to the one of a Denjoy counter-example.

\end{enumerate}
}
  \newpage
\section{Introduction}
The  exact symplectic twist maps of the two-dimensional annulus\footnote{all these notions will be precisely defined is subsection \ref{sub21}} were studied for a long time because they represent (via a
symplectic change of coordinates) the dynamic of the generic symplectic diffeomorphisms of
surfaces near their elliptic periodic points (see \cite{Ch1}). One motivating  example of such  a map was introduced by Poincar\'e for the study of  the restricted 3-Body problem.\\

The   study of such maps was initiated by G.D.~Birkhoff in the 20's (see \cite{Bir1}).  Among other beautiful results, he proved the following one (see \cite{He2} too)~: 
\begin{theo} {\bf (G.D.~Birkhoff)}
Let $f$ be a symplectic twist map of the two-dimensional annulus. Then any essential curve that is invariant by $f$ is the graph of a Lipschitz map.
\end{theo}
In this statement, an {\em essential curve} is a simple loop that is not homotopic to a point. \\

Later, in the 50's, the K.A.M.~theorems provide the existence of some invariant curves  for sufficiently regular symplectic diffeomorphisms of surfaces near their elliptic fixed points (see \cite{Ko}, \cite{Arno}, \cite{Mo} and \cite{Ru}). These theorems provide also some essential invariant curves for the symplectic twist maps that are close to the completely integrable ones. These K.A.M. curves are all very regular (at least $C^3$, see \cite{He2}). \\

But general invariant curves for general symplectic twist maps have no reason to be so regular. The example of the simple pendulum (see \cite{Ch2}) shows us that an invariant curve can be non-differentiable at one point: the separatrix of the simple pendulum has an angle at the hyperbolic fixed point. In \cite{He2} and \cite{Arna1}, some other examples are given of symplectic twist maps that have an non-differentiable essential invariant curve that contains some periodic points.\\

In all these examples, the non-differentiability appears at the periodic points. A natural question is then:\\
{\bf Question.} {\em Does a symplectic twist map  exist that has an essential invariant curve that is non-differentiable at a non-periodic point?}\\
A related question is the following one, due to J.~Mather in \cite{EKMY}:\\
{\bf Question. (J.~Mather)} {\em Does there exist an example of a symplectic $C^r$ twist map with an essential invariant curve that is not $C^1$ and that contains no periodic point (separate question for each $r\in [1, \infty]\cup\{ \omega\}$)?}\\

Let us point out that such an invariant essential curve  cannot be too irregular~:
\begin{enumerate}
\item[$\bullet$]  firstly, Birkhoff theorem implies that this curve has to be the graph of a Lipschitz map; hence, by Rademacher theorem, it has to be differentiable above a set that has full Lebesgue measure;
\item[$\bullet$] secondly, I proved in \cite{Arna1} that this curve has too be $C^1$ above a $G_\delta$ subset of $\T$ that has full Lebesque measure\footnote{The precise definition of $C^1$ in this context will be given in subsection \ref{sub21}}.
\end{enumerate}
We will prove:
\begin{thm}\label{T1}
There exists a symplectic $C^1$ twist map $f$ of the annulus that has an essential invariant curve $\Gamma$ such that:
\begin{enumerate}
\item[$\bullet$]  $\Gamma$ contains no periodic points;
\item [$\bullet$] the restriction $f_{|\Gamma}$ is $C^0$-conjugated to a Denjoy counter-example;
\item[$\bullet$] if $\gamma~: \T\rightarrow \R$ is the map whose $\Gamma$ is the graph, then $\gamma$ is $C^1$ at every point except along the projection of one orbit, along which $\gamma$ has distinct right and left derivatives.
\end{enumerate}
\end{thm}

This lets open Mather's question for $r\geq 2$ and also the following question:\\
{\bf Question.} {\em Does a symplectic twist map  exist that has an essential invariant curve that is non-differentiable and that is such that the dynamic restricted to this curve is minimal?}\\

Before giving the guideline of the proof, let us comment on some related results. \\

Finding invariant curves with no periodic points that are less regular that the considered symplectic twist map has been a challenging problem for a long time. \\
First at all, as an application of K.A.M. theorems, for a fixed diophantine rotation number and a 1-parameter smooth family (for example the standard family) such that the invariant curve disappears, it is classical that the ``last invariant curve'' is not $C^\infty$, even if the dynamic is $C^\infty$~: if this happens, by using K.A.M. theorems the curve cannot disappear\dots\\
Secondly, M.~Herman built in \cite{He2} some $C^2$ symplectic twist maps that have a $C^1$-invariant curve on which the dynamic is conjugated to the one of a Denjoy counter-example; such a curve cannot be $C^2$. In \cite{EKMY}, J.~Mather asks if such a $C^3$ counter-example exists.\\

To build our counter-example, we will use a family of symplectic twist maps that was introduced by M.~Herman in \cite{He2}. These maps are defined  by~:
$$f_\varphi: \T\times\R\rightarrow \T\times\R; (\theta, r)\mapsto (\theta+r, r+\varphi(\theta +r)).$$
where $\varphi~: \T\rightarrow \R$ is a $C^1$ map such that $\int_{\T}\varphi(\theta)d\theta=0$.\\
As noticed by M.~Herman, the main advantage of this map  is the following one.  We denote a lift of $g:\T \rightarrow \T$   by $\tilde g~:\R\rightarrow\R$. Then the graph of $\psi~: \T\rightarrow \R$ is invariant by $f_\varphi$ if and only if we have:
\begin{enumerate}
\item[$\bullet$] $g=Id_{\T}+\psi$ is an orientation preserving homeomorphism of $\T$;
\item[$\bullet$] $Id_\R +\frac{1}{2} \varphi=\frac{1}{2}\left( \tilde g+ \tilde g ^{-1}\right)$.
\end{enumerate}
Hence, in order to answer to Mather's question, we just have to find $g=Id_\T+\psi: \T\rightarrow \T$ that is an increasing non-differentiable homeomorphism of $\T$ with no  periodic points such that $\tilde g +\tilde g^{-1}$ is $C^1$.\\
We begin by using a $C^1$ Denjoy counter-example (see \cite{He1} and \cite{He2} for precise constructions) $g=Id_\T+\psi~:\T\rightarrow \T$. The non-wandering set of $g$ is then a Cantor subset $K$ of $\T$ such that $g_{|K}$ is minimal. We then consider a point $x_0\in \T\backslash K$ that is not in the Cantor subset $K$ and its orbit $(x_k)_{k\in\Z}=(g^k(x_0))_{k\in\Z}$. Then we modify $g$ in a neighborhood of this orbit   in such a way that the new homeomorphism $h~: \T\rightarrow\T$ coincides with $g$ along the orbit of $x_0$ and  is $C^1$ at every point but the orbit of $x_0$. At every point $x_k$ of the orbit of $x_0$, we assume that $h$ has some left and right derivatives, denoted by $\beta_k^l$ and $\beta_k^r$ such that: $\forall k\in\Z, \beta_k^r+\frac{1}{\beta_{k-1}^r}= \beta_k^\ell+\frac{1}{\beta_{k-1}^\ell}$.\\
If now we define $\varphi~: \T\rightarrow \R$ by~: $Id_\R +\frac{1}{2} \varphi=\frac{1}{2}\left( \tilde h+ \tilde h ^{-1}\right)$, then $\varphi$ is $C^1$ at every point of $\T$ but the orbit of $x_0$, has a right and left derivative along the orbit of $x_0$ and verifies (we denote by $\varphi_r'$ and $\varphi_\ell'$ the right and left derivative of $\varphi$): $\forall k\in \Z, \varphi'_r(x_k)= \beta_k^r+\frac{1}{\beta_{k-1}^r}-2= \beta_k^\ell+\frac{1}{\beta_{k-1}^\ell}-2=\varphi'_\ell(x_k)$. Hence $\varphi$ is differentiable. Roughly speaking, the left and right derivatives of $h$ at $x_k$ and $x_{k-1}$ are balanced in the formula that gives $\varphi$. This idea that the irregularities of $h$ and $h^{-1}$ are balanced in the formula that gives $\varphi$ was the one that used M.~Herman in \cite{He2} to construct his Denjoy counterexample for a $C^2$ symplectic twist map. 

If we choose carefully $h$, we will see  that $\varphi$ is in fact $C^1$. 
Let us point out that things are not as simple as they seem to be, and the choice of the sequences $(\beta_k^\ell)$ and $(\beta_k^r)$ is a delicate process as we will explain in the next section.

\section{Proof of theorem \ref{T1}} 
\subsection{Generalities about twist maps and other topics} \label{sub21}
\begin{nota}
\noindent $\bullet$ $\T=\R/\Z$ is the circle.

\noindent $\bullet$ $\A=\T\times \R$ is the annulus and an element of $\A$ is denoted
by $(\theta, r)$.

\noindent $\bullet$ $\A$ is endowed with its usual symplectic form, $\omega=d\theta\wedge dr$ and its usual Riemannian metric.

\noindent  $\bullet$ $\pi: \T \times \R \rightarrow\T$ is the first projection and $\tilde\pi:
\R^2\rightarrow \R$ its lift. 

\noindent $\bullet$ if $\alpha\in\T$, $R_\alpha: \T\rightarrow \T$ is the rotation defined by $R_\alpha (\theta)=\theta+\alpha$.

\end{nota}

\begin{defin} A $C^1$ diffeomorphism $f: \A\rightarrow \A$ of the annulus that is isotopic
to identity  is a {\em positive twist map} (resp. {\em negative twist map}) if, for any given lift $\tilde f: \R^2\rightarrow
\R^2$ and for every
$\tilde\theta\in\R$, the maps $r\mapsto \tilde\pi\circ \tilde f(\tilde\theta,r)$ is an increasing (resp decreasing)   diffeomorphisms. A {\em
twist map} may be positive or negative. 
\end{defin}
Then the maps $f_\varphi$ that we defined at the end of the introduction are positive symplectic twist maps.

\begin{defin} Let $\gamma~:\T\rightarrow\R$ be a continuous map. We say that $\gamma$ is $C^1$ at $\theta\in \T$ is there exists a number $\gamma'(\theta)\in\R$ such that, for every sequences $(\theta^1_n)$ and $(\theta_n^2)$ of points of $\T$ that converge to $\theta$ such that $\theta_n^1\not=\theta_n^2$, then:
$$\lim_{n\rightarrow \infty} \frac{\gamma(\theta_n^1)-\gamma(\theta_n^2)}{\theta^1_n-\theta^2_n}=\gamma'(\theta) $$
where we denote by $\theta_n^1-\theta_n^2$ the unique number that represents $\theta_n^1-\theta_n^2$ and that belongs to $]-\frac{1}{2}, \frac{1}{2}]$.
\end{defin}
If we assume that $\gamma$ is differentiable at every point of $\T$, then this notion of $C^1$ coincides with the usual one (the derivative is continuous at the considered point).
\subsection{Denjoy counter-example}
Following \cite{He2} p. 94, we define a Denjoy counter-example in the following way.\\
We assume that $\alpha\notin\Q/\Z$ , $\delta>0$ and that $C>>1$. Then we introduce:
$$\ell_k=\frac{a_C }{(|k|+C)(\log(|k|+C))^{1+\delta}}$$
where $a_C $ is chosen such that $\displaystyle{ \sum_{k\in\Z}\ell_k=1   }$. We use a $C^\infty$ function  $\eta~:\R\rightarrow \R$ such that $\eta\geq 0$, ${\rm support}(\eta)\subset [\frac{1}{4}, \frac{3}{4}]$ and $\int_0^1\eta(t)dt=1$. We define $\eta_k$ by~: $\eta_k(t)=\eta\left( \frac{t}{\ell_k}\right)$. Then we have: $\int_0^1\eta_k(t)dt=\ell_k$. Moreover, there exist two constants $C_1$, $C_2$, that depend  only on $\eta$, such that~:
$$C_1\leq |\eta_k|\leq C_2;\quad \frac{C_1}{\ell_k}\leq |\eta_k'|\leq \frac{C_2}{\ell_k}.$$
We assume now that $C>>1$ is great enough so that:
$$\forall k\in\Z, \left| \frac{\ell_{k+1}}{\ell_k}-1\right|C_2<1.$$
Then the map $g_k~: [0, \ell_k]\rightarrow [0, \ell_{k+1}]$ defined by $g_k(x)=\int_0^x \left(1+\left(\frac{\ell_{k+1}}{\ell_{k}}-1\right) \eta_k(t)\right)dt$ is a $C^\infty$ diffeomorphism such that $g_k(\ell_k)=\ell_{k+1}$.

There exists a Cantor subset $K\subset \T$ that has zero Lebesgue measure and that is such that the connected components of $\T\backslash K$, denoted by $(I_k)_{k\in\Z}$, are on $\T$ in the same order as the sequence $(k\alpha)$ and such that ${\rm length}(I_k)=\ell_k$.\\
Let us recall what is the semi-conjugation $j:\T\rightarrow \T$ of the Denjoy counter-example to the  rotation $R_\alpha$. Il $x\in \{ k\alpha; k\in\Z\}$, then we define~: $j^{-1}(x)=\int_0^xd\mu(t)$ where $\mu$ is the probability measure $\displaystyle{\mu =\sum_{k\in \Z}\ell_k\delta_{k\alpha}}$, $\delta_{k\alpha}$ being the Dirac mass at $k\alpha$. Then $j:\T\rightarrow\T$ is a continuous map with degree 1 that preserves the order on $\T$ and that is such that $j(I_k)=k\alpha$.\\
Then there is  a $C^1$ diffeomorphism $g:\T\rightarrow \T$ that fix $K$, is such that $K$ is the unique minimal subset for $g$, has for rotation number $\rho(g)=\alpha$, verifies $j\circ g=R_\alpha\circ j$. If $k\in\Z$, we introduce the notation: $g_{|I_k}=g_k$; then we have: $g_k(I_k)=I_{k+1}$. Following \cite{He2} again, we can assume that: $g_k'=g'_{|I_k}=\left(1+\left(\frac{\ell_{k+1}}{\ell_{k}}-1\right) \eta_k\right)\circ R_{-\lambda_k}$ where $R_{\lambda_k}(I_k)=[0, \ell_k]$  and that $g_k: I_k\rightarrow I_{k+1}$ is defined by~: $g_k=R_{\lambda_{k+1}}\circ h_k\circ R_{-\lambda_k}$.\\
Let us point out two facts that will be useful: $\displaystyle{\lim_{|k|\rightarrow \infty}\| g'_k-1\|=0
}$ and: \\
$\forall \theta\in K, g'(\theta)=1$.
\subsection{Modification of the Denjoy counter-example $g$}
We choose $x_0\in I_0$ and we consider its orbit $(x_k)_{k\in\Z}=(g^k(x_0))_{k\in\Z}$. \\
Then we will build a perturbation $h$ of $g$ such that $g_{|K}=h_{|K}$ and $\forall k\in\Z, h(x_k)=g(x_k)=x_{k+1}$.

 \begin{nota}
\begin{enumerate}
\item[$\bullet$] for all $k\in \Z$, we have: $I_k=]a_k, b_k[$, $L_k=]a_k, x_k]$ and $R_k=[x_k, b_k[$;
\item[$\bullet$] $\chi: \R\rightarrow \R$ is defined by: $\chi=\tilde g +\tilde g^{-1}-2Id_\R$;
\item[$\bullet$] for all $k\in\Z$, we denote: $\alpha_k=g'(x_k)$ and $m_k=2+\chi'(x_k)$.
\end{enumerate}
\end{nota}
Because of the definition of $\chi$, we have then: $\forall k\in\Z, \alpha_k+\frac{1}{\alpha_{k-1}}=m_k$.  

\begin{nota}
For every parameter $m\in\R$, let $\Phi_m: ]0, +\infty[\rightarrow ]-\infty, m[$ be defined by~: $\Phi_m(t)=m-\frac{1}{t}$.
\end{nota}
Let us notice that every function $\Phi_m$ is an increasing diffeomorphism.  Moreover:
\begin{enumerate}
\item[$\bullet$] if $m<2$, then $\Phi_m$ has no fixed points and~: $\forall t, \Phi_m(t)<t$ and $\displaystyle{\lim_{n\rightarrow +\infty}\Phi_m^n(t)=-\infty}$;
\item[$\bullet$] if $m=2$, $1$ is the only fixed point of $\Phi_m$. Moreover: if $t>1$, then $1<\Phi_m(t)<t$ and $\displaystyle{\lim_{n\rightarrow +\infty}\Phi_m^n(t)=1}$; if $t<1$, then $\Phi_m(t)<t$ and $\displaystyle{\lim_{n\rightarrow +\infty}\Phi_m^n(t)=-\infty}$;
\item[$\bullet$] if $m>2$, $\Phi_m$ has two fixed points, $p_-<p_+$; if $t>p_+$, then $p_+<\Phi_m(t)<t$ and 
$\displaystyle{\lim_{n\rightarrow +\infty} \Phi_m^n(t)=p_+}$; if $p_-<t<p_+$, then $p_-<t<\Phi_m(t)<p_+$ and 
$\displaystyle{\lim_{n\rightarrow +\infty} \Phi_m^n(t)=p_+}$; if $t<p_-$, then $\Phi_m(t)<t$ and $\displaystyle{\lim_{n\rightarrow +\infty}\Phi_m^n(t)=-\infty}$.
\end{enumerate}
  We have: $\forall k\in \Z, \alpha_k=\Phi_{m_k}(\alpha_{k-1})$.  \\
  
  Let us now choose $\beta_0^L>\alpha_0$ and $\beta_0^R>\alpha_0$ (each of them is then denoted by $\beta_0$). As every $\Phi_m$ is increasing and defined on $]0,+ \infty[$, we can define $(\beta_n)_{n\geq 0}$ in the following way: $\beta_{n+1}=\Phi_{m_n}(\beta_n)$. Then $\forall n\geq 0, \beta_n>\alpha_n>0$.
  \begin{lemma}\label{L1}
  We have: $\displaystyle{\lim_{n\rightarrow +\infty} \beta_n=1}$.
  \end{lemma}
  \demo Let us recall that: $\displaystyle{\lim_{n\rightarrow +\infty} \alpha_n=1}$. We deduce that $\displaystyle{\liminf_{n\rightarrow+\infty} \beta_n\geq 1}$ and  that $\displaystyle{\lim_{n\rightarrow +\infty} m_n=2}$.\\
   Let us fix $\varepsilon>0$; then there exists $N>0$ such that: $\forall n\geq N, m_n\leq 2+\varepsilon$. Then, for all $n\geq N$, we have: $\beta_{n+1}=\Phi_{m_{n+1}}(\beta_n)=\Phi_{2+\varepsilon}(\beta_n)-(2+\varepsilon-m_{n+1})<\Phi_{2+\varepsilon}(\beta_n)$. Using the fact that $\Phi_{2+\varepsilon}$ is increasing, we easily deduce: $\forall n\geq 0, \beta_{N+n}\leq \Phi_{2+\varepsilon}^n(\beta_N)$.  We know that $(\Phi_{2+\varepsilon}^n(\beta_N))_{n\in\N}$ has a limit, and because $\displaystyle{\liminf_{n\rightarrow+\infty} \beta_n\geq 1}$ this limit cannot be smaller than $1$. Hence
 $\displaystyle{\lim_{n\rightarrow +\infty} \Phi_{2+\varepsilon}^n(\beta_n)=p_+(\varepsilon)}$ if we denote the greatest fixed point of $\Phi_{2+\varepsilon}$ by $p_+(\varepsilon)$. We deduce that $\displaystyle{\limsup_{n\rightarrow +\infty}\beta_n\leq p_+(\varepsilon)}$. We have: $p_+(\varepsilon)=\frac{2+\varepsilon+\sqrt{\varepsilon(4+\varepsilon)}}{2}$, hence: $\displaystyle{\lim_{\varepsilon\rightarrow 0^+}p^+(\varepsilon)=1}$.\\
   Finally, we have proved that $\displaystyle{\limsup_{n\rightarrow +\infty}\beta_n\leq 1}$. As   $\displaystyle{\liminf_{n\rightarrow +\infty}\beta_n\geq 1}$, we deduce the lemma.
   
  \enddemo
  In a similar way, we choose $0<\beta_{-1}^L<\alpha_{-1}$ and $0<\beta_{-1}^R<\alpha_{-1}$ and we denote each of them by $\beta_{-1}$. As every $\Phi_m^{-1}$ is increasing and defined on $]-\infty, m[$, we can define $(\beta_k)_{k\leq -1}$ by $\beta_{k-1}=\left(\Phi_{m_k}\right)^{-1}(\beta_k)$. Then : $\forall k\leq -1, \beta_k\leq \alpha_k$.
\begin{lemma}
  We have: $\displaystyle{\lim_{k\rightarrow -\infty} \beta_k=1}$.
\end{lemma}
The proof of this lemma is similar to the one of lemma \ref{L1}.
Now, we can assume:
\begin{enumerate}
\item[$\bullet$] $\beta_0^R\not=\beta_0^L$;
\item[$\bullet$] $\beta_0^L+\frac{1}{\beta_{-1}^L}=\beta_0^R+\frac{1}{\beta_{-1}^R}$.
\end{enumerate}
Then we denote this last quantity by $\tilde m_{0}=\beta_0^L+\frac{1}{\beta_{-1}^L}=\beta_0^R+\frac{1}{\beta_{-1}^R}$. Because $\beta_0>\alpha_0$ and $\beta_{-1}<\alpha_{-1}$, we necessarily have: $\tilde m_0>m_0$.

Now $(\beta_k^L)$ (resp. $(\beta_k^R))$ is a good candidate to be the left (resp. right) derivative of $h$ along the orbit $(x_k)$.
\begin{prop}\label{P1}
There exists an orientation preserving homeomorphism $h:\T\rightarrow \T$ such that:
\begin{enumerate}
\item[$\bullet$]  $h_{|K}=g_{|K}$ and $\forall k\in\Z, h(x_k)=g(x_k)$;
\item[$\bullet$]  $h$ and $h^{-1}$ are $C^1$ at every point  of $\T$ but the orbit of $x_0$;
\item[$\bullet$] $h$ and $h^{-1}$ have some right and left derivative at every point $x_k$ of the orbit of $x_0$, and $h'_R(x_k)=\beta_k^R$, $h'_L(x_k)=\beta_k^L$. Moreover, $h_{|R_k}$ and $h_{|L_k}$ are $C^1$.
\end{enumerate}
\end{prop}
\begin{remk}
If   proposition \ref{P1} is true,  then theorem \ref{T1} is proved: if $h=Id_\T+\psi$ and $\varphi=\tilde g+\tilde g^{-1}-2Id_\R$, then the graph of $\psi$ is invariant by $f_\varphi$ and the dynamic of $f_\varphi$ restricted to this graph is the one of a Denjoy counter-example with $\alpha$ as rotation number. Moreover, $\psi$ is non-differentiable along the orbit of $x_0$ but $\varphi$ is $C^1$.  Indeed, as $g$ and $g^{-1}$ are, $\varphi$ is $C^1$ at every point of $\T$ but the orbit of $x_0$. Moreover, the restriction of $\varphi$ to each interval $L_k=]a_k, x_k]$ or $R_k=[x_k, b_k[$ is $C^1$. To prove that $\varphi$ is $C^1$, we then just have to prove that the right and left derivatives are equal along the orbit of $x_0$. We have:
\begin{enumerate}
\item[$\bullet$] if $k\not=0$, $\varphi'_L(x_k)=\beta_k^L+\frac{1}{\beta_{k-1}^L}-2=m_k-2=\beta_k^R+\frac{1}{\beta_{k-1}^R}-2=\varphi'_R(x_k)$;  
\item[$\bullet$] if $k=0$, $\varphi'_L(x_0)=\beta_0^L+\frac{1}{\beta_{-1}^L}-2=\tilde m_0-2=\beta_0^R+\frac{1}{\beta_{-1}^R}-2=\varphi'_R(x_0)$.
\end{enumerate}
Hence $\varphi$ is $C^1$.
\end{remk}
Let us now prove proposition \ref{P1}. We modify $g$, or rather its derivative, in each interval $L_k$ and $R_k$ in the following way. Let us notice that:
 $\displaystyle{\lim_{|k|\rightarrow +\infty} |g'_{|[a_k, b_k]}-1|=0}$;
  $\displaystyle{\lim_{|k|\rightarrow +\infty}\beta_k^L=\lim_{|k|\rightarrow +\infty}\beta_k^R=1}$;
  $g'_{|K}=1$.\\
Then on each interval $L_k=]a_k, x_k]$, we replace $g'$ by a continuous function $\delta_k:]a_k, x_k]\rightarrow \R_+^*$  such that:
\begin{enumerate}
\item\label{pt1} $\delta_k(x_k)=\beta_k^L$;
\item\label{pt2} $\delta_k$ coincide with $g'$ in a neighborhood of $a_k$;
\item\label{pt3} $\int_{a_k}^{x_k}\delta_k=\int_{a_k}^{x_k}g'=\tilde g(x_k)-\tilde g(a_k)$;
\item\label{pt4} $\forall t\in L_k, |\delta_k(t)-1|\leq \max\{ |g_{|L_k}-1|, |\beta_k^L-1|\}+\frac{1}{1+|k|}$.
\end{enumerate}
To build $\delta_k$, we just have to replace $g'$  between $x_k-\varepsilon_k$ and $x_k$ by some affine function and then to modify slightly $g'$ elsewhere in $L_k$ to rectify the value of the integral. If $\varepsilon_k$ is small enough, than the change in the integral is very small and we have the last inequality (but of course the slope of the affine function can be very great, so the perturbation of $g$ that we build in not small in $C^2$ topology).\\
We then define $h_{|L_k}$ by: 
$$\forall t\in [a_k, x_k], h(t)= g(a_k)+\int_{a_k}^t\delta_k (s)ds.$$
We proceed in similar way to define $h_{|R_k}$ and we obtain similar properties.
Moreover, we ask: $h_{|K}=g_{|K}$. 

Then $h$ is continuous. By construction, its restriction to every interval $I_k$ is continuous. Moreover, we have $\displaystyle{\lim_{|k|\rightarrow +\infty} |h_{|I_k}-g(a_k)|=0}$ (because of point \ref{pt4} and the fact that $\displaystyle{\lim_{|k|\rightarrow \infty}{\rm length}(I_k)=0}$). We deduce that $h$ is continuous at every point of $K$. Moreover, $h$ is orientation preserving and injective by construction. Hence $h$ is an orientation preserving homeomorphism of $\T$.

Moreover, $h_{|K}=g_{|K}$ by construction and $\forall k\in\Z, h(x_k)=g(x_k)$ by point \ref{pt3}.

By construction, $h_{|R_k}$ and $h_{|L_k}$ are $C^1$ (and then the same is true for $h^{-1}$), $h'_R(x_k)=\beta_k^R$ and $h'_L(x_k)=\beta_k^L$.

Let us now prove that $h$ is differentiable at every point of $K$ and that $h'_{|K}=1$. We consider $y\in K$ and a sequence $(y_n)$ that converge to $y$ and that is such that: $\forall n, y_n\not=y$.  We want to prove that $\displaystyle{\lim_{n\rightarrow +\infty} \frac{h(y_n)-h(y)}{y_n-y}=1}$. Considering eventually different cases, we can assume that $(y_n)$ tends to $y$ from above. Then there are two cases:
\begin{enumerate}
\item[$\bullet$] either $y=a_k$ for some $k$. Then we have the conclusion by point \ref{pt2};
 \item[$\bullet$]  or $y$ is accumulated from above by a sequence $(a_{k_n})_{n\in\N}$ that are left ends of intervals $I_{j_k}$.
\end{enumerate}
Because $g$ is $C^1$ and   $g'_{|K}=1$, for every $\varepsilon>0$ there exists $\eta>0$ such that for every $z\in [y, y+\eta[$, then $|g'(z)-1|<\varepsilon$. Let us assume that $n$ is big enough such that $y_n\in[y, y+\eta[$. There are three cases:
\begin{enumerate}
\item[$\bullet$]  $y_n\in K$. Then there exists $z\in [y, y+\eta[$ such that: $$ \left| \frac{h(y_n)-h(y)}{y_n-y}-1\right|= \left|\frac{g(y_n)-g(y)}{y_n-y}-1\right|=\left|g'(z)-1\right|<\varepsilon;$$
\item[$\bullet$]  $y_n\in L_k$ for some $k$. Then there exist $z, z'\in [y, y+\eta[$ such that:
$$  \frac{h(y_n)-h(y)}{y_n-y} =  \frac{h(y_n)-h(a_k)+g(a_k)-g(y)}{y_n-y} =\frac{y_n-a_k}{y_n-y}h'(z)+\frac{a_k-y}{y_n-y}g'(z').$$
\item[$\bullet$]  $y_n\in R_k$ for some $k$. Then there exist $z, z'\in [y, y+\eta[$ such that:
$$  \frac{h(y_n)-h(y)}{y_n-y} =  \frac{h(y_n)-h(x_k)+g(x_k)-g(y)}{y_n-y} =\frac{y_n-x_k}{y_n-y}h'(z)+\frac{x_k-y}{y_n-y}g'(z').$$
Because $g$ is $C^1$ and because of point \ref{pt4}, if $\eta$ is small enough, then $h'(z)$ and $g'(z')$ are close enough to 1, and their  barycentre is close to 1 two.

 \end{enumerate}
 Hence we have prove that $h$ is derivable along $K$ and that $h'_{|K}=1$.

Because $g$ is $C^1$ and $g'_{|K}=1$, because of point \ref{pt4}, $h$ is $C^1$ on $K$ and $h'_{|K}=1$.

Then $h$ satisfies all the conclusions of proposition \ref{P1}.

\end{document}